\documentclass[aos,MSNbibl,nameyear,dvips]{arximspdf}
\usepackage{dcolumn}

\doi{10.1214/14-AOS1175F} 
\referstodoi{10.1214/13-AOS1175}
\volume{42}
\issue{2}
\pubyear{2014}
\firstpage{509}
\lastpage{517}

\makeatletter
\newcommand{\df}{\mathrm{df}}
\newcolumntype{d}[1]{D{.}{.}{#1}}
\newcommand{\Exp}{\operatorname{Exp}}
\newcommand{\eqref}[1]{(\ref{#1})}
\newcommand{\0}{\mathbf{0}}
\newcommand{\y}{\mathbf{y}}
\newcommand{\I}{\mathbf{I}}
\newcommand{\X}{\mathbf{X}}
\newcommand{\N}{\mathcal{N}}
\newcommand{\Beps}{\bolds{\varepsilon}}
\newcommand{\sigmahat}{\hat{\sigma}}
\makeatother

\begin{document}
\begin{frontmatter}

\title{Discussion: ``A significance test for the lasso''}
\runtitle{Discussion}

\begin{aug}
\author{\fnms{A.} \snm{Buja}\corref{}\ead[label=e1]{buja.at.wharton@gmail.com}}
\and
\author{\fnms{L.} \snm{Brown}\ead[label=e2]{lbrown@wharton.upenn.edu}}
\runauthor{A. Buja and L. Brown}
\affiliation{University of Pennsylvania}
\address{Statistics Department\\
The Wharton School\\
University of Pennsylvania\\
471 Jon M. Huntsman Hall\\
Philadelphia, Pennsylvania 19104-6340\\
USA\\
\printead{e1}\\
\phantom{E-mail:\ }\printead*{e2}} 
\pdftitle{Discussion of ``A significance test for the lasso''}
\end{aug}

\received{\smonth{1} \syear{2014}}



\end{frontmatter}

It is rare in our field that we can speak of a true discovery, but
this is one such occasion. It is an unexpected result that the steps
by which variables enter a lasso path permit a basic statistical test
with a simple null distribution that is asymptotically valid and has
good finite-sample properties. This test may become standard, and
maybe it should simply be called ``\emph{the} lasso test'' because it
is difficult to conceive of a form of inference more intimately tied
to the lasso.\footnote{Moreover, the term ``covariance test'' is
misleading because it is not covariance that is being tested.}

The authors use forward stepwise variable selection as their straw man,
and this for good reason because the $t$-tests on which stepwise
selection builds are essentially a heuristic abuse of the testing
framework that entirely ignores the effects of selection. The lasso
tests, by contrast, account for selection and shrinkage that is
implicit in the lasso. Insight into how this is possible is one of
the many byproducts of this innovative and thought provoking article.

One of the beauties of the authors' article is Lemma 3 where it is
shown that the transformation of sorted test statistics $|z^{(j)}|$
using the formula $|z^{(j)}| (|z^{(j)}| - |z^{(j+1)}|)$ generates
quantities that have limiting null distributions $\Exp(1/j)$ as \mbox{$p
\rightarrow\infty$}. By comparison, $|z^{(j)}|$ keeps growing under the
null at the rate $\sqrt{2 \log p}$. It is a quite remarkable fact
that for this particular series of transformed statistics a limiting
distribution can be obtained under competition by an unlimited number
of null predictors, that is, $p \rightarrow\infty$.

\subsection*{From straw man to competitor: Forward stepwise}

The authors' Section~2.2 convincingly documents that naive $t$-tests
are fatally flawed when used in the standard forward stepwise
selection routine to test the conditional null hypothesis that the
current selection contains all nonzero coefficients. In the example
of their Figure~1, the authors consider testing the first selected
predictor among ten orthogonal predictors.\vadjust{\goodbreak} Assuming $\sigma$ known,
the $t$-statistic becomes a $z$-statistic, and the null distribution
of $z^2$ for any of the predictors is $\chi_1^2$. If, however, the
tested predictor has been chosen to maximize explanatory\vspace*{1pt} power, then
its proper null distribution is not $\chi_1^2$ but
$ \max_{j=1,\ldots,10} \chi_{1(j)}^2$, where $\chi_{1(j)}^2$ are ten
independent copies of $\chi_1^2$. The authors' Figure~1(a)
illustrates the obvious fact that this distribution is stochastically
\emph{much} larger than the naive null distribution~$\chi_1^2$.

Once this is recognized, however, there are various ways to account
for the effects of selection in the forward stepwise procedure. In
what follows, we briefly outline how forward stepwise selection can be
provided with inference that is conditionally valid given the
selection path taken thus far, just like the lasso test, but unlike
the lasso test, the inference is guaranteed to be strictly valid for
finite samples and also for arbitrary collinearities. In detail,
consider selection stage $k$ where the set $A$ of selected predictors
has size $|A| = k - 1$. The number of remaining predictors is
$p - (k - 1)$, and we denote these by $\X_{k},\ldots,\X_{p}$. In
view of the fact that the predictors in $A$ have already been
included, we need versions of the remaining predictors that are
adjusted for the predictors in $A$, and we denote these versions by
$\X_{j \cdot A}$ ($j = k,\ldots,p$). Unlike the authors, we will not
assume that $\sigma$ is known but that it must be estimated by some
$\sigmahat= \sigmahat(\y)$ with $\df_{\mathrm{err}}$ degrees of freedom
(usually from the RSS of the full model) and which we can assume to be
stochastically independent of all $\langle\X_j, \y\rangle$.
Importantly, we use this single $\sigmahat$ for all $t$- and
$F$-statistics and never recompute it from any submodel. This is
important in order to enable simultaneous inference to solve the
multiplicity problem of selection [\citet{Beretal13}, Section~4.1]. In
order to test the strongest among the remaining $p - (k - 1)$
predictors under the null hypothesis that $A$ contains all predictors
with true nonzero slopes, one can proceed in one of the following
ways:
\begin{itemize}
\item \textit{Exact inference based on max-$|t|$}: Assuming that the selected
predictor at each step is the one with the most significant
$t$-statistic if added to the model $A$, the appropriate test
statistic is
\begin{equation}
\label{eqmax-t} t_{\max}(\y):= \max_{j=k,\ldots,p}
\bigl|t^{(j)}(\y)\bigr|\qquad\mbox{where } t^{(j)}(\y):=
\frac{\langle\X_{j \cdot A}, \y\rangle}{\|\X_{j \cdot A}\| \sigmahat
(\y)}.
\end{equation}
The null distribution of $t_{\max}(\y)$ under the assumption that
all remaining predictors have zero slopes can be approximated by
simulating $t_{\max}(\Beps)$ for $\Beps\sim\N(\0,\I_n)$, while for
small numbers of remaining predictors there exists software to
perform numerical integration. The correct $p$-value is
$\mathbf{P}[t_{\max}(\Beps) > t_{\max}(\y)]$. This is the brute-force
approach that correctly accounts for any finite sample size and
arbitrary collinearities. It is only weakness is that it assumes
homoskedastic normal errors whose variance is properly estimated by
$\sigmahat^2$; first-order correctness of the full model does not
need to be assumed if such a $\sigmahat$ is available [\citet{Beretal13}, Sections~2.2~and~3].
\item \textit{Bonferroni correction to naive inference}: Use naively the
$t_{\df_{\mathrm{err}}}$-distribution for $t_{\max}(\y)$, but adjust the
significance level by dividing it by $p - (k - 1)$, or else
adjust the naive $p$-value by multiplying it by $p - (k - 1)$.
This approach is conservative but provides excellent approximations
for nearly orthogonal predictors.
\item \textit{Scheff\'e simultaneous inference}: The Scheff\'e method can be
used to provide simultaneous inference for all linear combinations
of the remaining coefficients, which trivially includes all of the
remaining coefficients. Scheff\'e-adjusted \mbox{$p$-}values are obtained by
treating $t_{\max}^2(\y)/(p - (k - 1))$ as distributed
according to $F_{p-(k-1),\df_{\mathrm{err}}}$. This approach is obviously too
conservative but it is easy to obtain as an alternative when the
Bonferroni correction fails due to strong collinearity.

\item \textit{$F$-tests of remaining variation}: This method is not strictly a
test of selected predictors but it has a touch of the obvious in
that a significant test result suggests that more predictors should
be included. The method consists of performing an $F$-test of each
submodel $A$ within the full model. Note that for orthogonal
predictors the $F$-statistic at stage $k$ is $\sum_{j=k,\ldots,p}
(t^{(j)} )^2 / (p - (k - 1))$, where $t^{(j)} =
t^{(j)}(\y)$ is defined in \eqref{eqmax-t}. The null distribution
is again $F_{p-(k-1),\df_{\mathrm{err}}}$. One can give the method an
interpretation in the spirit of Scheff\'e simultaneous inference: A
significant $F$-test at a given stage means that there exists a
linear combination of the remaining coefficients that is
statistically significant. This, then, suggests continuing with
inclusion of another term. One stops stepwise inclusion when the
$F$-test indicates that there does not exist a linear combination of
the remaining predictors that accounts for significant variation in
the response.
\item \textit{Lemma}~2 \textit{tests for stepwise}: Lemma 2 can be used for stepwise
selection, but it, too, is not strictly a test of selected
predictors because it depends not only on the strongest but the
second strongest remaining predictor as well. The test statistic is
$t_{\max}\cdot(t_{\max}-t_{\max-1})$, where $t_{\max-1}$ stands for
the second largest in magnitude among $t$-statistics of remaining
predictors. According to Lemma~2, for orthogonal predictors this
test statistic has an approximate $F_{2,\df_{\mathrm{err}}}$-distribution. As
in the case of the $F$-test method, a statistically significant
outcome of the Lemma~2 test indicates that more predictors are
needed. The power implications of this choice of test statistic
are not clear at this point, although the authors provide some
tentative simulation results in their Figure~4, which in their
example seems to indicate no drastic differences in power between
the Lemma~2 statistic and the $t_{\max}$ statistic.
\end{itemize}
All types of $p$-values for the sequence of forward stepwise inclusions
are shown in Table~\ref{tab1} for the full wine quality data (the
authors show their results for a half sample, hence some disagreements
with our results). The exact method is based on 99,999 null\vadjust{\goodbreak}
replicates. Computation of the whole table took just eight seconds in
spite of the simulations at each step for the exact method.

\begin{table}
\def\arraystretch{0.9}
\tabcolsep=0pt
\caption{Comparative results for the various conditional $p$-values in
forward stepwise selection applied to the full wine quality data}\label{tab1}
\begin{tabular*}{\tablewidth}{@{\extracolsep{\fill}}@{}lcd{2.4}cccccc@{}}
\hline
&&& \multicolumn{6}{c@{}}{$\bolds{p}$\textbf{-values}}\\[-6pt]
&&& \multicolumn{6}{c@{}}{\hrulefill}
\\
\textbf{Step} & \textbf{Predictor} & \multicolumn{1}{c}{\textbf{$\bolds{t}$-stats}} & \textbf{Naive} & \textbf{Exact} & \textbf{Bonfer} & \textbf{Scheffe} & \textbf{$\bolds{F}$-tests} &\textbf{Lemma 2}\\
\hline
\phantom{0}1 & Alcohol             &  23.7216 & 0.0000 & 0.0000 & 0.0000 & 0.0000 & 0.0000 & 0.0000\\
\phantom{0}2 & Volatile\_acidity     & 14.9676 & 0.0000 & 0.0000 & 0.0000 & 0.0000 & 0.0000 & 0.0000\\
\phantom{0}3 & Sulphates            &  6.8479 & 0.0000 & 0.0000 & 0.0000 & 0.0000 & 0.0000 & 0.0000\\
\phantom{0}4 & Total\_sulfur\_dioxide &  4.4237 & 0.0000 & 0.0001 & 0.0001 & 0.0125 & 0.0000 & 0.4136\\
\phantom{0}5 & Chlorides            &  4.3749 & 0.0000 & 0.0001 & 0.0001 & 0.0080 & 0.0000 & 0.0011\\
\phantom{0}6 & pH                   &  3.7544 & 0.0002 & 0.0011 & 0.0011 & 0.0291 & 0.0010 & 0.0062\\
\phantom{0}7 & Free\_sulfur\_dioxide  &  2.3878 & 0.0171 & 0.0726 & 0.0853 & 0.3369 & 0.1370 & 0.0915\\
\phantom{0}8 & Citric\_acid          &  1.0633 & 0.2878 & 0.6540 & 1.0000 & 0.8893 & 0.6124 & 0.6309\\
\phantom{0}9 & Residual\_sugar       &  0.7818 & 0.4344 & 0.7528 & 1.0000 & 0.8938 & 0.6705 & 0.8429\\
10 &Fixed\_acidity        &  0.5071 & 0.6122 & 0.7829 & 1.0000 & 0.8794 & 0.6250 & 0.8190\\
11 &Density              &  0.8266 & 0.4086 & 0.4066 & 0.4086 & 0.4086 & 0.4086 & \\
\hline
\end{tabular*}   \vspace*{-3pt}
\end{table}

In conclusion, forward stepwise selection can be richly endowed with
valid statistical inference. It does not deserve to be seen as the
poor ``step child'' of the lasso.\vspace*{-2pt}

\subsection*{Issues with the application of lasso tests}

If the history of the $t$-test is a guide, the lasso test will give us
some quirks and curiosities to ponder. Part of the historic learning
curve in connection with $t$-tests was the experience that
occasionally two predictors can be both statistically insignificant
when they appear jointly in a model, but when one of them is removed
the other is boosted to statistical significance (one of the joys of
collinearity). As a consequence, it was understood that it is not a
good idea to simultaneously remove all insignificant predictors from a
model. This in turn led to the invention of stepwise selection
procedures, then to the lasso, and now to the article at hand.

As an example of an issue to ponder about the lasso test, there is the
notion of a random null hypothesis. In our own work on valid
post-selection inference [\citet{Beretal13}], we faced a similar issue
and referee questions: what does it mean to provide valid inference in
a random model? This question is unavoidable when the models in which
tests are to be performed are the result of a random selection process
such as a stepwise, all-subsets or lasso variable selection procedure.
The way the issue was resolved in our work was by providing protection
for all possible null hypotheses that \emph{could} have been selected,
hence the selection procedure provides only a lens to randomly
focus
on one of many null hypotheses whose validity of inference has been
insured beforehand. This is not so for lasso tests: they are truly
conditional starting with the second selection of a predictor.\vadjust{\goodbreak}

An issue arising from sequential conditionality can be illustrated by
scanning some of the data examples provided by the authors: for the
prostate cancer data in Table~1, we will say that after (conditional
on) including the first variable (\texttt{lcavol}) we have evidence at a
level just barely missing significance 0.05 that the second variable
(\texttt{lweight}) carries signal. In the subsequent four steps, the
added variables do not provide evidence that they carry signal given
the inclusion of the respective previous variables. But, conditional
on including six predictors, the seventh (\texttt{lcp}) gives evidence
again just barely missing significance 0.05 that it, too, carries
signal. This could, of course, be a false rejection, but if Lemma~3 and
Theorem~1 of the article are a guide, the sequence of null
distributions becomes tighter [$\Exp(1)$, $\Exp(1/2)$, $\Exp(1/3),\ldots,$ for
orthogonal predictors] under repeated null inclusions. As a
consequence, if we assume that the insignificant $p$-values in steps~3
through 6 correspond to true nulls, then at step 7 the true null
distribution might be as tight as $\Exp(1/5)$, meaning that the value of
the covariance statistic [which is not shown but we figure to be about
3.134 if based on $F(2,58)$] could be multiplied by a factor up to 5,
resulting in a $p$-value as low a 0.0000036 as opposed to 0.051. The
conclusion is that there is something wrong with the assumption that
steps 3 through 6 are null inclusions in spite of their
insignificances. This might be something to chew on.\looseness=-1

The effect just described does not seem to be isolated as it appears
again, in milder form, in the authors' training half-sample of the
wine quality data (Table~5): if we assume there that steps 4 and 5 are
null inclusions, then step 6 with a \mbox{$p$-}value of 0.076 could have an
effective test statistic larger by a factor up to 3, amounting to an
effective $p$-value as low as 0.00044. Again, the conclusion is that
the assumption that steps 4 and 5 were null inclusions is wrong in
spite of their insignificances. In our replication with the full wine
quality data, the Lemma~2
test features one erratic jump into insignificance at step 4 before
resuming with significance for two more steps thereafter.

The erratic and somewhat trend-less behavior of sequences of
conditional \mbox{$p$-}values down a lasso path is an issue with which
practitioners will struggle. It would be desirable to smooth the
sequences so they show more a trend than erratic jumps. One potential
approach to this problem could be some form of bootstrap smoothing or
bagging. Here is an attempt: we adopted a crude criterion using the
0.05 threshold to chose as estimated model size the largest step
number whose lasso $p$-value sequence up to that point remains below
0.05. Shown in Table~\ref{tab2} are the cumulative counts of model
sizes; for example, for model size 4, 991 out of 1000 bootstrap
resamples generated a lasso test sequence whose first four $p$-values
remained below 0.05, hence the estimated model size is 991 out of 1000
times at least~4.

\begin{table}
\def\arraystretch{0.9}
\tabcolsep=0pt
\caption{$p$-values summaries for bootstrap based on 1000 resamples
applied to the full wine quality data}\label{tab2}
\begin{tabular*}{\tablewidth}{@{\extracolsep{\fill}}@{}lcccccccc@{}}
\hline
\textbf{Step} &\textbf{1} &\textbf{2} &\textbf{3} &\textbf{4} &\textbf{5} &\textbf{6} &\textbf{7}& \textbf{8}\\
\hline
Cumul\# $p$-values $< 0.05$ &  1000  & 1000   & 1000   & 991    & 936    & 818    & 576    & 350 \\
Median $p$-value           & 0.000  & 0.000  & 0.000  & 0.008  & 0.057  & 0.048  & 0.169  & 0.370\\
\hline
\end{tabular*}
\end{table}

While we are somewhat dubious regarding the meaning of these numbers,
they do seem to suggest that lasso $p$-values are erratic. For example,
the large observed $p$-value of 0.537 at step 5 (total\_sulfur\_dioxide)
may have been a fluke because 936 bootstrap resamples produce a
$p$-value below 0.05 up to and including step~5, and the median $p$-value
at step 5 is 0.057. Even at step 4, the observed $p$-value of 0.173
seems excessive in view of the fact that the median $p$-value at that
step is 0.008. In summary, there is evidence that the observed
$p$-values require some kind of processing because they do seem to
behave somewhat unpredictably. We have not tried the methods proposed
by Grazier G'Sell et al. (\citeyear{GG13}) which form aggregates of the observed
$p$-values to control FDR. Aggregation may help somewhat, but intuition
suggests that aggregating the observed series of $p$-values may not
achieve sufficient smoothing; one may have to shake up the data
repeatedly and aggregate the results to achieve greater stability of
inferential conclusions.

To be fair to the authors, they do not actually make recommendations
how to use the lasso $p$-values. Theirs is a technical article that
lays out the theory and concepts, but it does not propose a
methodology. This, however, is done in a companion article by Grazier G'Sell et al. (\citeyear{GG13}), a must-read for
anyone who cares about actually using lasso tests. While this is a
wide-ranging article, its Section~5 discusses sequential selection
rules for ``FDR control for the lasso in \emph{nonidealized
settings}'' (emphasis added by us). Perusing this material
shortened the present discussion considerably because the authors have
already worked through many of the issues that arise when lasso tests
meet practice, some of which we were about to raise on our own. For
efficient dissemination of the news we allow ourselves to quote some
striking insights (spoiler alert) and add our own comments:
\begin{itemize}
\item``\emph{Breakdown of the} $\Exp(1/l)$ \textit{behavior}$\ldots.$ \textit{In finite
samples}, \textit{the} $\Exp(1/l)$ \textit{behavior becomes unreliable for larger}
$l$, \textit{leading the corresponding statistics to be larger than
expected}.'' In response to this issue, the authors limit the
look-ahead in one of their rules to $l \le$ 5 or 10. Our earlier
observations on the prostate cancer and wine quality data involved
$l = 5$ and $l = 3$, respectively, and are therefore within the
authors' limits. So here the solution is simple: do not expect the
$\Exp(1/l)$ behavior to hold for long stretches.
\item``\emph{Intermingling of signal and noise variables$\ldots.$ The
hypotheses made by Lockhart et al.} (\citeyear{Lo13}) \textit{prevent this from
happening asymptotically}, \textit{but the assumptions can still break down
in practice.}'' Asymptotic theory assumes that the
\mbox{signal}\vadjust{\goodbreak}
variables end up in the active set before the noise variables, which
with sufficient data will be the case with high probability but will
be doubtful in any given data situation. There is apparently no
good solution to this problem as the authors report that this can
render their rule to be anti-conservative. A general qualm we have
with present day's excitement over sparseness is that in our
experience data are rarely sparse in the sense that signal variables
stick out from a background of noise variables like a mesa (as in
Lockhart et al.'s Theorem 1 where a $\sigma\sqrt{2 \log p}$
threshold is assumed). Signal tends to peter out gradually and will
be sparse only in the sense that for large $p$ only a small fraction
of predictors have signal that is detectable. As a result, we find
ourselves in need of making trade-offs that will always be
unsatisfactory to some when deciding where to come down on the scale
from conservative to liberal. The fact that, in practice, signal
tends to exist at all scales does not invalidate the use of tests
for zero signal, but users of sequences of lasso tests (or stepwise
tests) have to contend with the fact that the sequence will
ultimately hit a gray zone where signal and noise variables start to
mingle.
\item``\emph{Correlation in $X\ldots.$ the null distributions of Lockhart
et al.} (\citeyear{Lo13}) \textit{begin to break down when $X$ has high correlation.}''
The take-home message is that, in practice, collinearity cannot be
defined away. We may need some diagnostics to help us decide what
form and degree of collinearity invalidates the $\Exp(1)$ form of the
null distribution. We may also have to live with trade-offs again,
as when the true mean response is $X_1+X_2$ but $X_1$ and $X_2$ are
highly collinear in relation to the sample size $n$ and the noise
level $\sigma$, in which case model selection procedures will make a
random choice between the two predictors.
\item``\emph{The appropriateness of FDR as an error criterion becomes
questionable when $X$ is highly correlated. If a noise variable is
highly correlated with a signal variable}, \textit{should we consider it to
be a false selection}? \textit{This is a broad question that is beyond the
scope of this paper}, \textit{but is worth considering when discussing
selection errors in problems with highly correlated $X$.}'' This
comment speaks to us like no other. Our work on valid
post-selection inference [\citet{Beretal13}] is FWER-based, but we
wondered what a FDR-based version would look like. We, too, decided
that FDR does not even make sense for similar reasons: if some
predictors form a cluster (are mutually highly collinear), then
there are many ways of making a selection error that really amounts
to the same error, whereas for a predictor that is nearly orthogonal
to all others there is only one way to make this selection error.
As a consequence, counting selection errors and forming rates does
not seem meaningful in the presence of collinearity; the FDR concept
needs adjusting, but it is not obvious how.
\end{itemize}

\subsection*{Larger issues in statistical inference}

Finally, we wish to step back and discuss some larger issues. While
the authors' article is a tremendous advance, it is a first and
necessary step on a long path to solve larger problems:
\begin{longlist}[(1)]
\item[(1)] Lasso tests do assume an underlying linear model with Gaussian
errors. At some point, we may need tests that do not require this
assumption. Generally speaking, we will need statistical inference
that is valid under model misspecification [\citet{Bujetal}]. The
dangers from misspecification should increase as data with $p > n$
become common place because diagnosing nonlinearity and
heteroskedasticity will become impossible, yet their effects on
sampling variability and inference will persist, just better
concealed due to their undiagnosability in the $p > n$ regime.
\item[(2)] Each class of tests, be they lasso tests or stepwise tests, can
ultimately be augmented in such a way that they control FWER or FDR
if used sequentially [Grazier G'Sell et al. (\citeyear{GG13})]. However, this
assumes for their validity that data analysts obey a protocol
whereby they commit a priori to one and only one selection method,
lasso, for example, and nothing else. Now consider the more
realistic situation in which a data analyst tries both, lasso and
forward stepwise selection, and decides based on gut feeling or
informal devices such as plots which of the two to use: if the data
analyst is honest at heart and clear in his mind, he will realize
that he faces a meta-selection problem. Compounding the problem is
that he may not even have followed a generalizable rule in his
decision in favor of lasso or stepwise. How are we to evaluate such
practice and its effects on statistical inference? One of the
benefits of our approach in \citet{Beretal13} is that it sets
analysts free to do experimenting with selection methods to their
hearts' content, followed by meta-selection according to any rule or
none---subsequent inference will still be valid.
\item[(3)] The preceding point opens up the bigger issue of informal
methods, often graphical, that are used for exploratory data
analysis and model diagnostics. Such methods often inform data
analysts in fruitful ways to guide them to more meaningful analyses,
but they may have insidious effects on subsequent inference.
Analysts may have no feelings of dishonesty and may not be aware
that they are biasing the analysis and modeling process in
unaccountable ways. It just seems like the reasonable thing to do
to prevent nonsense from happening. We have tried to introduce a
small measure of inference in the EDA and model diagnostics process
in \citet{Bujetal09} and \citet{Wicetal10}, but the larger
question remains unanswered: what is the compounded effect of the
many informal activities at all stages of data analysis on
statistical inference?
\item[(4)] Empirical research has taken to statistics with a vengeance in
less than half a century. Yet, empirical research suffers from a
systemic malady that is well reflected by \citet{Ioa05} piece
with the provocative yet realistic title ``Why most published
research findings are false.'' The culprit of first order is most
likely publication bias, also called the ``file drawer problem,''
that is, the fact that negative results tend not to see the light of
publication. A culprit of second order we hypothesize to be the
fact of unaccounted data analytic activity, ranging from
meta-selection among variable selection methods to the use of
informal EDA and \mbox{diagnostics} methods. It may just be the case that
the most expert and thorough data analysts are also the ones who
produce the most spurious findings in applied statistical work.
This should not be construed as a call to apply less competence and
abandon research into efficient statistical methods, but it should
be motivation to create statistical inference that integrates ever
more of the informal data analytic activities for which there is
currently no accounting. This is again some of the background of
our proposal in \citet{Beretal13} which provides valid
post-selection inference even if data analysts are arbitrarily
informal in their meta-selection of variable selection methods.
\end{longlist}

Returning to the occasion of this discussion, clearly the authors'
article represents an advance which, with suitable methodology, will
fill a large missing piece in statistical inference. We hope that
this and forthcoming pieces will ultimately coalesce into a larger
methodology that will account for data analytic activities which still
fall through the cracks of current best practice. We conclude by
thanking the authors for an inspiring article.



%

\printaddresses

\end{document}